\documentclass[ 12pt, leqno]{article}

\usepackage{amsmath,amssymb,amsthm}
\usepackage{tikz}
\usepackage{color}
\usepackage[toc]{appendix}
\usepackage{graphicx}
\usepackage{fancyhdr}
\usepackage{enumitem}
\usepackage{bbm}
\usepackage{parskip}
\usepackage{float}
\usepackage{chngpage}
\usepackage{calc}
\usepackage{bigints}
\usepackage{array}
\usepackage{booktabs}
\usepackage{rotating}
\usepackage{multirow}
\usepackage{adjustbox}
\usepackage{tabularx}
\usepackage{verbatim}
\usepackage{mathtools}
\usepackage{ragged2e}
\usepackage[makeroom]{cancel}
\usepackage{caption}
\usepackage{hyperref}
\usepackage{caption}
\usepackage{subcaption}
\usepackage{appendix}
\usepackage{pgfplots}
\pgfplotsset{compat=1.16}
\usepackage[english]{babel}
\usepackage{hyphenat}
\usepackage{accents}

\newtheorem{them}{Theorem}[section]

\newtheorem{lem}{Lemma}[section]

\newtheorem{remark}{Remark}[section]

\numberwithin{equation}{section}

\author{ Enrico Bernardi \thanks{Dipartimento di Scienze Statistiche Paolo Fortunati, Universit\`a di Bologna, Bologna, Italy. \textbf{e-mail}: enrico.bernardi@unibo.it} \and
Tatsuo Nishitani\thanks{Department of Mathematics, Graduate School of Science, Osaka University,
  Toyonaka, Osaka, Japan.  \textbf{e-mail}:nishitani@math.sci.osaka-u.ac.jp}
}

\date{****}

\title{Notes on tangent bicharacteristics and ill-posedness of the Cauchy problem}

\begin{document}

\maketitle

\def\dif{\partial}
\def\al{\alpha}
\def\be{\beta}
\def\ga{\gamma}
\def\te{\theta}
\def\bphi{{\bar\phi}}
\def\op#1{{\rm op}(#1)}

\def\lr#1{\langle{#1}\rangle}
\def\xig{\langle{\xi}\rangle_{\!\gamma}}
\def\R{{\mathbb R}}
\def\C{{\mathbb C}}
\def\N{{\mathbb N}}

\def\de{\delta}

\def\La{\Lambda}
\def\la{\lambda}
\def\ep{\epsilon}
\def\varep{\varepsilon}
\def\lxi{\langle{\xi'}\rangle}

\def\Ga{\Gamma}

\begin{abstract}
We exhibit a family of second-order hyperbolic differential operators
presenting spectral transition of the Hamilton map. As a consequence
we prove that the Cauchy problem is not locally solvable at the origin
in  Gevrey classes of order greater than some fixed value. The main
feature of these operators is that they may all have bicharacteristics
tangent to the double manifold.
\end{abstract}

%%%%%%%%%%
%%%%%%%%%%%%%%%
%

%%%%%%%%%%%
%%%%%%%%%%%
\section{Introduction}
Ever since in \cite{Ivrii-Fact} a remarkable correlation was
discovered between  the structure of the well-posedness
of the Cauchy problem for hyperbolic operators with double
characteristics and the
behavior of the Hamilton flow of their principal symbols, a
number of papers have investigated this phenomenon, see e.g. \cite{NCE} ,
\cite{BN1} and particularly \cite{Nibook}, Chapters 6,7 and 8. It has thus become
clear that besides the standard Ivrii-Petkov-H\"ormander conditions
bearing on the lower order terms of the operator, this new feature
on how the null bicharacteristics of the principal symbol fail or
happen  to be
attracted by the set of double points plays a major part in the final
assessment of when the $ C^{\infty} $ well-posedness of the Cauchy
problem for the operator can be ascertained.
Up to now however  only cases
where the spectral type of the Hamilton map of the principal symbol on
the double manifold is constant have been analyzed.

It is therefore the purpose
of this paper to examine a family of examples extending the known
ill-posedness results to symbols admitting a change of spectral
behavior, a so called transition of type, along the manifold of double
points.
Filtering the class of $ C^{\infty} $ functions via the Gevrey scale,
it then becomes possible to pinpoint with a certain precision when, in
terms of the Gevrey index, the failure of well-posedness takes place. As
proven below in Theorem \ref{thm:ill-posed} this index will partly depend 
also on
the geometric behavior of those bicharacteristics landing tangentially
on the double manifold. The examination of these model families, all having
analytic coefficients and null subprincipal symbols, relies on
the geometrical analysis performed in \cite{BN:arxiv} where the
symplectic structure of the transition was thoroughly studied. Furthermore a major technical ingredient turns out be a careful consideration of the
regions where the Stokes coefficients attached to a family of ordinary
differential equations in the complex domain may vanish.

%In this note we prove
%
\begin{them}
\label{thm:ill-posed} Consider
\[
P=-D_0^2+2x_1D_0D_n+D_1^2+x^2_1(x_1+\varep x_{r}^l)D_n^2,\quad 2\leq r\leq n-1
\]
in $|x_i|\ll 1$ where $l\in\N$ and $\varep=\pm 1$. There exists $N>1$, independent of $l$,  such that the Cauchy problem for $P$ is not locally solvable at the origin in the Gevrey class of order 
\begin{equation}
\label{eq:sisu}
s>\max{\{5, 1+N/l\}}
\end{equation}
it is even more so for $C^{\infty}$.
\end{them}
\begin{remark}
\label{rem:b}\rm Note that the subprincipal symbol of $P$ is identically $0$. Actually, no matter how  the $b_j\in\C$ are chosen, the same conclusion holds for $P+\sum_{j=0}^nb_jD_j$.
\end{remark}
Let us briefly explain the motivation for considering this family of operators.
To that end, consider a slightly more general family of operators;
\[
\tilde P=-D_0^2+2x_1D_0D_n+D_1^2+x_1^2(x_1^k+\varep x_r^l)D_n^2,\quad |x_i|\ll 1,\;\; 2\leq r\leq n-1
\]
where $k\in \N$ and $P$ in Theorem \ref{thm:ill-posed} corresponds to the case $k=1$. When $k=1$ and $\varep=0$, that is
% {\color{blue} It is clear that a bicharacteristic $\gamma(s)=(x(s),\xi(s))$ with $x_p(s)=0$ of
%
\[
P_{mod}=-D_0^2+2x_1D_0D_n+D_1^2+x_1^3D_n^2
\]
which is a model of non-effectively hyperbolic of type 2 {\it without spectral transition} %that is,
%
%\[
%W(\rho)={\rm Ker}F_{p}^2(\rho)\cap{\rm Im} F_p^2(\rho)\neq \{0\},\quad \rho\in \Sigma
%\]
%
but a {\it bicharacteristic tangent to the doubly characteristic set}, which is given explicitly by
\begin{equation}
\label{eq:bicha}
\begin{cases}
x_1=-x_0^2/4,\quad x_n=x_0^5/80,\quad \xi_1=x_0^3/8,\quad \xi_n=1\\
x_j=\text{constant},\quad \xi_j=\text{constant},\quad 2\leq j\leq n-1
\end{cases}
\end{equation}
parametrized by $x_0$. A detailed study is found in \cite{BN1}, \cite{Nibook} where it was shown that 
the Cauchy problem for $P_{mod}$ is not locally solvable at the origin in the Gevrey class $s>5$ (a part of arguments in \cite[Section 5.2]{BN1} should read as \cite[Section 6.2]{Nibook}) and that for more general operator enjoying the same properties as $P_{mod}$ the Cauchy problem is well-posed at the origin in the Gevrey class of order less than $5$ under the Levi condition. Note that if $k=1$ for any $l$, $\varep$ there exists a bicharacteristic tangent to the doubly characteristic set of 
\[
\tilde p=-\xi_0^2+2x_1\xi_0\xi_n+\xi_1^2+x_1^2(x_1+\varep x_r^l)\xi_n^2,\quad \tilde P=\op{\tilde p}
\]
because $\tilde p$ is independent of $\xi_r$, setting $x_r=0$ in \eqref{eq:bicha} we obtain such a bicharacteristic unless $l=1$, if $l=1$ the $\xi_r$ is obtained by integrating the corresponding equation.  

Since {\it a spectral transition occurs} for $\tilde p$ we briefly describe it. Making a change of coordinates
\[
y_j=x_j,\quad 0\leq j\leq n-1, \quad y_n=x_n+x_0x_1
\]
then in the new local coordinates we have
\[
\tilde P=\op{-\xi_0^2+(\xi_1+x_0\xi_n)^2+x_1^2(1+x^k_1+\varep x_r^l)\xi_n^2},\quad |x_i|\ll1.
\]
The doubly characteristic set $\Sigma$ near  $\xi_n\neq0$ is given by
\[
\Sigma\cap\{\xi_n\neq 0\}=\{\xi_0=\xi_1+x_0\xi_n=0, x_1=0\}.
\]
Denoting 
\begin{gather*}
\varphi_1=-x_1(1+x_1^k)\xi_n/(1+(1+k/2)x_1^k)=-x_1\al_k(x_1)\xi_n,\quad \al_k(0)=1 
\end{gather*}
one can write
\[
\tilde p=-(\xi_0+\varphi_1)(\xi_0-\varphi_1)+(\theta_k+\varep x_r^l/\al_k^2)\varphi_1^2+\varphi_2^2,\quad \varphi_2=\xi_1+x_0\xi_n
\]
where 
\begin{equation}
\label{eq:rei:3:b}
\theta_k=(1+k)x_1^k+k^2x_1^{2k}(1+x_1^k)^{-1}/4
\end{equation}
which is of normal form (see \cite[Section 2]{BN:arxiv}). 
Since $(\theta_k+\varep x_r^l/\al_k^2)|_{\Sigma}=\varep x_r^l$ it follows from \cite[Proposition 2.1]{BN:arxiv} that the spectral transition depends only on $\varep x_r^l$ which is independent of $k\in\N$, that is $\tilde p$ is effectively hyperbolic where $\varep x_r^l<0$ and non-effectively 
hyperbolic of type 1 or type 2  where  $\varep x_r^l>0$ or $\varep x_r^l=0$.
In particular, when $\varep=1$ and $l$ is even, as mentioned above, $\tilde p$ is noneffectively hyperbolic of type 1 in $\Sigma\cap\{x_r\neq 0\}$ and of type 2 on $\Sigma\cap \{x_r=0\}$ such that the spectral transition occurs on $\Sigma\cap\{x_r=0\}$ which is independent of $k\in\N$. On the other hand, when $k=1$ there is a bicharacteristic tangent to $\Sigma\cap\{x_r=0\}$ as remarked above and the Cauchy problem for $\tilde P$ is not locally solvable at the origin in the Gevrey class of order greater than \eqref{eq:sisu} in view of Theorem \ref{thm:ill-posed}, while if $k\geq 2$ thanks to \cite[Section 5]{BN:arxiv} this $\tilde p$ admits an elementary factorization which implies that there is no tangent bicharacteristic and that the Cauchy problem for $\tilde P$ is $C^{\infty}$ well-posed at the origin.

%
%%%%%%%%%%%%%%
\section{Proof of Therem}

Since  $2\leq r\leq n-1$ has no special meaning, without loss of
generality we may fix $r=n-1$ in what follows. 

%%%%%%%%%%%%%%%%%
\subsection{Sibuya solutions}

Let
\[
P_{\rho}=-D_0^2+2x_1D_0D_n+D_1^2+x_1^2(x_1+\varep \rho^{-(N+2)}x_{n-1}^l)D_n^2
\]
where $\rho>0$ is a positive large parameter and $N>2$ will be determined later and we look for a solution to $P_{\rho}U=0$ of the form
\begin{equation}
\label{eq:teigi:U}
U(x;\al,\rho)=\exp{\big(-i\rho^5x_n-\frac{i}{2}\zeta\rho x_0\big)}w(\rho^2x_1)
\end{equation}
where $x=(x_0,x_1,x_n)$ and $\zeta\in \C$. Note that $P_{\rho}U=0$ gives
\begin{gather*}
-w''(\rho^2x_1)+\big\{\zeta x_1\rho^2 +x_1^2(x_1+\varep \rho^{-(N+2)}x_{n-1}^l)\rho^6-\zeta^2\rho^{-2}/4\big\}w(\rho^2x_1)\\
=-w''(y)+\big\{y^3+(\varep \rho^{-N}x_{n-1}^l)y^2+\zeta y-\zeta^2\rho^{-2}/4\big\}w(y)=0
\end{gather*}
with $y=\rho^2x_1$. In what follows, we regard $x_{n-1}\in\R$ as a small parameter. Denoting 
\begin{equation}
\label{eq:teigi:al}
\al=\varep \rho^{-N}x_{n-1}^l
\end{equation}
one can write
\begin{gather*}
y^3+\al y^2+\zeta y-\zeta^2\rho^{-2}/4=(y+\al/3)^3+(\zeta-\al^2/3)(y+\al/3)\\
-\al\zeta/3+2\al^3/27-\zeta^2\rho^{-2}/4
\end{gather*}
so that $w$ satisfies
\[
-w''(y-\al/3)+\big\{y^3+(\zeta-\al^2/3)y-\al\zeta/3+2\al^3/27-\zeta^2\rho^{-2}/4\big\}w(y-\al/3)=0.
\]
Consider
\begin{equation}
\label{eq:Sibu}
-Y''(x)+(x^3+\zeta x+\ep)Y(x)=0,\quad \zeta,\;\ep\in\C.
\end{equation}
This equation has a solution ${\mathcal Y}(x;\zeta,\ep)$ which is entire analytic in $(x,\zeta,\ep)$ and admits an asymptotic representation
\[
{\mathcal Y}(x;\zeta,\ep)\sim x^{-3/4}(1+\sum_{N=1}^{\infty}B_Nx^{-N/2})\exp{\{-E(x;\zeta)\}}
\]
uniformly on each compact set in the $(\zeta,\ep)$ space as $x$ goes to infinity in any closed subsector of the open sector $|\arg{x}|<3\pi/5$ where
\[
E(x,\zeta)=\frac{2}{5}x^{5/2}+\zeta x^{1/2}
\]
(see \cite[Chapter 2]{Sib}).
If we choose 
\begin{equation}
\label{eq:def:w}
w(y)={\mathcal Y}(y+\al/3;\zeta-\al^2/3, -\al\zeta/3+2\al^3/27-\zeta^2\rho^{-2}/4)
\end{equation}
then $U$ given by \eqref{eq:teigi:U} with \eqref{eq:def:w} satisfies $P_{\rho}U=0$.

Denote
\[
{\mathcal Y_k}(x;\zeta,\ep)={\mathcal Y_k}(\omega^{-k}x;\omega^{-2k}\zeta,\omega^{-3k}\ep),\;\; \omega=e^{i2\pi/5},\quad k=0,1,2,3,4,\quad {\mathcal Y_0}={\mathcal Y}
\]
then every ${\mathcal Y_k}(x;\zeta,\ep)$ solves \eqref{eq:Sibu} and tends to $0$ as $x$ tends to infinity along any direction in the sector
\[
\big|\arg{x}-2k\pi/5\big|<3\pi/5.
\]
We have
\[
{\mathcal Y_k}(x;\zeta,\ep)=C_k(\zeta,\ep){\mathcal Y_{k+1}}(x;\zeta,\ep)-\omega {\mathcal Y_{k+2}}(x;\zeta,\ep)
\]
(see \cite[Chapter 5]{Sib}). Recall that if $C_0(\zeta,\ep)=0$ then ${\mathcal Y}(x;\zeta,\ep)={\mathcal Y_0}(x;\zeta,\ep)$ behaves 
\begin{equation}
\label{eq:Y:x<0}
\begin{split}
|{\mathcal Y}(x;\zeta,\ep)|\approx \exp{(-2x^{5/2}/5)},\quad \R\ni x\to +\infty,\\
|{\mathcal Y}(x;\zeta,\ep)|\approx \exp{({\mathsf{Im}\,\zeta}|x|^{1/2})},\quad \R\ni x\to -\infty
\end{split}
\end{equation}
if $(\zeta,\ep)$ remains to be bounded.
Look for a $\zeta=\zeta(\rho^{-1}, x_{n-1})$ satisfying
\[
C_0(\zeta-\al^2/3, -\al\zeta/3+2\al^3/27-\zeta^2\rho^{-2}/4)=0,\quad {\mathsf{Im}}\big(\zeta-\al^2/3)<0.
\]
First note that
\begin{equation}
\label{eq:rouche}
\begin{split}
|C_0(\zeta-\al^3/3,-\al\zeta/3+2\al^3/27-\zeta^2\rho^{-2}/4)-C_0(\zeta,-\zeta^2\rho^{-2}/4)|\\
\leq C_1|\al|^3+C_2|\al\zeta/3-2\al^3/27|
\leq C'|\al|
\end{split}
\end{equation}
for bounded $|\zeta|$ and $|x_{n-1}|$. Consider $C_0(\zeta,-\zeta^2s/4)=0$ and recall
\begin{lem}{\rm(\cite[Proposition 6.1]{Nibook})}
\label{lem:ze:0}
$C_0(\zeta,0)$ has a zero $\zeta_0$ with negative imaginary part.
\end{lem}

Let $\zeta_0$ be in Lemma \ref{lem:ze:0}. Since $C_0(\zeta,0)$ is not identically zero (see \cite[Chapter 5]{Sib}), thanks to the Weierstrass preparation theorem one can write
\[
C_0(\zeta,-\zeta^2s/4)=c(\zeta,s)\big((\zeta-\zeta_0)^m+\sum_{j=1}^ma_{k}(s)(\zeta-\zeta_0)^{m-k}\big)=c(\zeta,s)Q(\zeta,s)
\]
in a neighborhood of $(\zeta_0,0)$ with some $m\in\N$ where $c(\zeta_0, 0)\neq 0$ and $a_j(0)=0$. Recall that zeros of $Q(\zeta,s)$ are given by convergent Puiseux series 
\[
\zeta_j(s)=\sum_{k\geq 0}\zeta_{jk}(s^{1/p_j})^k,\quad 1\leq j\leq m, \quad \zeta_{j0}=\zeta_0,\quad p_j\in\N,\quad |s|<\delta.
\]
Choose one of them, say $\zeta_1(s)$ then it is clear that one can find $c>0$ and $\nu_1>0$ such that
\[
|\zeta_j(s)-\zeta_1(s)|\geq c|s|^{\nu_1}
\]
unless $\zeta_j(s)$ is identically equal to $\zeta_1(s)$ which proves that there exist $\nu>0$, $c_i>0$ and $\mu>0$ such that
\[
|C_0(\zeta, -\zeta^2s/4)|\geq c_2|s|^{\mu}\quad\text{if}\quad |\zeta-\zeta_1(s)|=c_1|s|^{\nu}
\]
where necessarily $\mu>\nu$. 
From now on we denote $\eta(\rho^{-1})=\zeta_1(\rho^{-1})$.
Since
\begin{gather*}
|\al|\leq |x_{n-1}|^l\rho^{-N}\leq C\rho^{-N},\quad  |x_{n-1}|\leq c'
\end{gather*}
taking \eqref{eq:rouche} into account and choosing $N$ such that $N>\mu$ one can apply Rouch\'e's theorem to find $\zeta=\zeta(\rho^{-1}, x_{n-1})$ satisfying
\begin{equation}
\label{eq:w:to:Y}
\begin{split}
C_0(\zeta-\al^3/3,-\al\zeta/3+2\al^3/27-\zeta^2\rho^{-2}/4)=0,\\
 |\zeta(\rho^{-1}, x_{n-1})-\eta(\rho^{-1})|<c_1\rho^{-\nu}.
 \end{split}
\end{equation}
Since $
|\zeta(\rho^{-1},x_{n-1})-\zeta_0|\leq |\zeta(\rho^{-1}, x_{n-1})-\eta(\rho^{-1})|+|\eta(\rho^{-1})-\zeta_0|$ when $\rho\to\infty$  $\zeta(\rho^{-1},x_{n-1})$ converges to $\zeta_0$ uniformly in $|x_{n-1}|<c'$.
Therefore from \eqref{eq:Y:x<0} one can find $C>0, c>0$ such that
\begin{equation}
\label{eq:cal:Y}
\begin{split}
\big|{\mathcal Y}(x;\zeta-\al^2/3,-\al\zeta/3+2\al^3/27-\zeta^2\rho^{-2}/4)\big|\leq Ce^{-c|x|^{1/2}},\\ 
\zeta=\zeta(\rho^{-1}, x_{n-1}),\quad \al=\varep\rho^{-N}x_{n-1}^l
\end{split}
\end{equation}
holds on $x\in\R$  uniformly in large $\rho$ because ${\mathsf{Im}}\big(\zeta(\rho^{-1},x_{n-1})-\al^3/3\big)\leq -c_1$ with $c_1>0$ and $\al$, $\zeta$ remains in a bounded set for large $\rho$. 
%

%%%%%%%%%%%%%
\subsection{Proof of Theorem}
Consider the Cauchy problem
\begin{equation}
\label{eq:Gev:tekisetu}
\begin{cases}
P_{\rho}v_{\rho}=0,\\
v_{\rho}(0, x')=0,\\
D_0v_{\rho}(0,x')=\phi_1(x_1)\tilde\phi(x'')\psi(x_{n-1})\theta(x_n)
\end{cases}
\end{equation}
with Gevrey Cauchy data $\phi_1, \psi, \theta\in \ga_0^{(s_1)}(\R)$ and $\tilde\phi\in \ga_0^{(s_2)}(\R^{n-3})$ where $x''=(x_2,\ldots,x_{n-2})$ and $s_1>s_2>1$ will be precised later. 
Note that the coefficients of $P_{\rho}$ are analytic (actually polynomials) in $x$, and whose radii of convergence and the modulus are independent of $\rho\geq 1$. Thus the existence domain of analytic solutions to the Cauchy problem (assured by the Cauchy-Kowalevsky Theorem) can be taken independent of polynomial Cauchy data and $\rho\geq1$. Then applying Holmgren's arguments  (see, for example, Mizohata \cite[Theorem 4.2]{Mizohata}) we conclude that there exists  
\begin{equation}
\label{eq:D:iki}
D_{\delta_0}=\{x\in \R^{1+n};|x'|^2+|x_0|<\delta_0\}
\end{equation}
independent of $\rho\geq1$ such that if $u\in C^2(D_{\delta})$ ($0<\delta\leq \delta_0$) satisifes $P_{\rho}u=0$ in $D_{\delta}$ and $\dif_{x_0}^ju(0,x')=0$ ($j=0,1$) on $D_{\delta}\cap \{x_0=0\}$ then $u=0$ in $D_{\delta}$. Moreover applying this uniqueness theorem, choosing the Cauchy data with small supports and taking $T>0$ small, for $v_{\rho}\in C^2(D_{\delta})$ (if exists) satisfying \eqref{eq:Gev:tekisetu} one can assume that 
\[
\{x;0\leq x_0\leq T, v_{\rho}(x)\neq 0\}\Subset K,\quad \rho\geq 1
\]
for any compact neighborhood $K$ of the origin given arbitrarily beforehand.

Here recall some facts from \cite[Section 6.3]{Nibook}. Let $h>0$ and $L$ be a compact set in $\R^n$. We say $f(x)\in \ga_0^{(s),h}(L)$ if $f(x)\in C_0^{\infty}(L)$ and
\begin{equation}
\label{eq:Gev:norm}
\sum_{\al}\sup_x h^{|\al|}|\dif_x^{\al}f(x)|/|\al|!^s
\end{equation}
is finite. Here $\ga_0^{(s),h}(L)$ is a Banach space equipped with the norm \eqref{eq:Gev:norm}. 
\begin{lem}
\label{lem:kyotuiki} {\rm (see, for example, \cite[Proposition 4.1]{Mizohata})}
Let $L$ be a compact set of $\R^n$ and $h>0$ be fixed. Assume that the Cauchy problem for $P=P_1$ is locally solvable in $\ga^{(s)}$ at the origin. Then there is $\delta>0$ such that for any $(u_0(x'), u_1(x'))\in (\ga_0^{(s),h}(L))^2$ there exists a unique $u(x)\in C^2(D_{\delta})$ verifying $Pu=0$ in $D_{\delta}$ and $D_0^ju(0,x')=u_j(x')$ on $D_{\delta}\cap \{x_0=0\}$.
\end{lem}
Replacing $L$ by a smaller one if necessary, we can assume that $L\subset D_{\delta}\cap \{x_0=0\}$. Since the linear map $(\ga_0^{(s),h}(L))^2\ni (u_0(x'), u_1(x'))\to u\in C^1(D_{\delta})$ is closed then from Banach's closed graph theorem it follows that this map is continuous. Let $K\subset D_{\delta}$ be a compact set such that $L\Subset K\cap\{x_0=0\}$. Then there is $C>0$ such that
\begin {equation}
\label{eq:futousiki}
|u|_{C^1(K)}\leq C\sum_{\al}\sup_{j,x'} h^{|\al|}|\dif_{x'}^{\al}u_j(x')|/|\al|!^s.
\end{equation}
Suppose that the Cauchy problem for $P_1$ is locally solvable in $\ga^{(s)}$ at the origin with $s$ verifying \eqref{eq:sisu}. Noting that $N/l<s-1$ and $s>5$ we choose $s>s_1>1$ and $s>s'_2>s_2>1$ so that
\begin{equation}
\label{eq:stos'}
s>s_1>5,\quad N/l<s-s'_2.
\end{equation}
Thanks to Lemma \ref{lem:kyotuiki} there is a solution $u_{\rho}\in C^2(D_{\delta})$ to the Cauchy problem
\begin{equation}
\label{eq:CP:rho=1}
\begin{cases}
P_1u_{\rho}=0,\\
u_{\rho}(0, x')=0,\\
D_0u_{\rho}(0,x')=\phi_1(x_1)\tilde\phi(x'')\psi(\rho^{N/l}x_{n-1})\theta(x_n).
\end{cases}
\end{equation}
If we denote $\tilde u=u_{\rho}(x_0,x_1,x'',\rho^{-N/l}x_{n-1},x_n)\in C^2(D_{\delta})$ it is clear that $\tilde u$ satisfies \eqref{eq:Gev:tekisetu} then from the uniqueness of the solution to the Cauchy problem \eqref{eq:Gev:tekisetu} holding in \eqref{eq:D:iki} for any $\rho\geq 1$ we have
\[
v_{\rho}(x)=u_{\rho}(x_0, x_1, x'', \rho^{-N/l}x_{n-1}, x_n).
\]
 Choose $\ep>0$ such that $h\ep<1$ and note that
\begin{gather*}
\frac{|\dif_{\tilde x}^{\al}\dif_{x_{n-1}}^kD_0u_{\rho}(0,x')|}{(k+|\al|)!^s}\leq \frac{CA^{k+|\al|}|\al|!^{s_1}k!^{s_2}\rho^{kN/l}}{(k+|\al|)!^s}\leq \frac{C'\ep^{k+|\al|}|\al|!^{s}k!^{s'_2}\rho^{kN/l}}{(k+|\al|)!^s}\\
\leq C'\ep^{k+|\al|}\frac{\rho^{kN/l}}{k!^{s-s'_2}}\leq C'\ep^{k+|\al|}e^{(s-s'_2)\rho^{N/l(s-s'_2)}}
\end{gather*}
where $\tilde x=(x_1,\ldots,x_{n-2}, x_n)$. Therefore we have
\begin{equation}
\label{eq:koko}
\sum_{k+|\al|}h^{k+|\al|}\sup \frac{|\dif_{\tilde x}^{\al}\dif_{x_{n-1}}^kD_0u_{\rho}(0,x')|}{(k+|\al|)!^s}\leq Ce^{(s-s'_2)\rho^{N/l(s-s'_2)}}.
\end{equation}
Hence with $\kappa=N/l(s-s'_2)<1$ we conclude from \eqref{eq:futousiki} that
\begin{equation}
\label{eq:Gev:katei:1}
\sup_{K}|u_{\rho}|+\sup_{K}|D_0u_{\rho}|\leq Ce^{\rho^{\kappa}}.
\end{equation}
Since $e^{i\zeta \rho T/2}$ is independent of $(x_0,x_1,x_n)$ considering $e^{i\zeta \rho T/2}U$ instead of $U$ we see that
\begin{gather*}
P_{\rho}U=0,\quad  U=e^{-i\rho^5x_n+\frac{i}{2}\zeta\rho (T-x_0)}w(\rho^2 x_1),\\
 w(y)={\mathcal Y}(y+\al/3;\zeta-\al^2/3,-\al\zeta/3+2\al^3/27-\zeta^2\rho^{-2}/4)
\end{gather*}
where $
\zeta=\zeta(\rho^{-1}, x_{n-1})$ and $ \al=\varep \rho^{-N}x_{n-1}^l$.
Note that we have
\begin{equation}
\label{eq:Gev:katei:b}
\begin{split}
\sup_{K}|v_{\rho}|\,(\leq \sup_{K}|u_{\rho}|), \:\;\sup_{K}|D_0v_{\rho}|\,(\leq \sup_{K}|D_0u_{\rho}|)\leq Ce^{\rho^{\kappa}}.
\end{split}
\end{equation}
We repeat the same argument as in \cite[Section 6.3]{Nibook}. 
\begin{gather*}
\int_0^T(P_{\rho}U, v_{\rho})dx_0=\int_0^T(U, P_{\rho}v_{\rho})dx_0+i(D_0U(T),  v_{\rho}(T))\\
+i(U(T), D_0v_{\rho}(T))-i(U(0), D_0v_{\rho}(0))-i(2x_1D_nU(T), v_{\rho}(T))
\end{gather*}
because $v_{\rho}(0)=0$. From this we have
\begin{equation}
\label{eq:green}
\begin{split}
(D_0U(T), v_{\rho}(T))+(U(T), D_0v_{\rho}(T))\\
-2(x_1D_n U(T), v_{\rho}(T))=(U(0), D_0v_{\rho}(0)).
\end{split}
\end{equation}
Here note that
\begin{gather*}
U(T)=e^{-i\rho^5x_n}w(\rho^2x_1),\quad 
 D_nU(T)=-\rho^5e^{-i\rho^5x_n}w(\rho^2x_1),\\
  D_0U(T)=-\rho \zeta(\rho^{-1},x_{n-1}) e^{-i\rho^5x_n}w(\rho^2x_1).
\end{gather*}
Then it follows from \eqref{eq:green} and \eqref{eq:Gev:katei:b} that there is $C>0$ such that
\[
|(U(0), D_0v_{\rho}(0))|\leq C\rho^5e^{\rho^{\kappa}}
\]
for $|w(\rho^2x_1)|$ is bounded on $\R$ uniformly in $\rho$ by \eqref{eq:cal:Y}.
On the other hand
\begin{gather*}
(U(0), D_0v_{\rho}(0))=\int_{\R^n} e^{-i\rho^5x_n+i\zeta(\rho^{-1}, x_{n-1})\rho T/2}w(\rho^2x_1)\\
\times \phi_1(x_1)\tilde\phi(x'') \psi(\rho^{N/l}x_{n-1})\theta(x_n)dx_1dx''dx_{n-1}dx_n\\
=\rho^{-N/l-2}\hat\theta(\rho^5)\int_{\R^{n-3}}\tilde\phi (x'')dx''\\
\times \int_{\R^2}e^{i\zeta(\rho^{-1}, x_{n-1})\rho T/2}w(x_1)\phi_1(\rho^{-2}x_1)\psi(x_{n-1})dx_1dx_{n-1}
\end{gather*}
where $\hat\theta$ is the Fourier transform of $\theta$. Consider
\[
I_{\rho}=\int_{\R^2}e^{i\zeta(\rho^{-1}, x_{n-1})\rho T/2}w(x_1)\phi_1(\rho^{-2}x_1)\psi(x_{n-1})dx_1dx_{n-1}.
\]
Write
\begin{gather*}
\int w(x_1)\phi_1(\rho^{-2}x_1)dx_1=\sum_{j=0}^2\frac{\rho^{-2j}}{j!}\phi_1^{(j)}(0)\int w(x_1)x_1^jdx_1+O(\rho^{-6}).
\end{gather*}
Thanks to \cite[Lemma 6.6]{Nibook} one can assume that there is $j$ ($0\leq j\leq 2$) such that
\[
\int {\mathcal Y}(x_1;\zeta_0,0)x_1^idx_1=0,\quad 0\leq i\leq j-1,\quad\int {\mathcal Y}(x_1;\zeta_0,0)x_1^jdx_1\neq 0.
\]
From \eqref{eq:w:to:Y} we have $
\big|\zeta(\rho^{-1}, x_{n-1})-\eta(\rho^{-1})\big|\rho\leq c_1\rho^{1-\nu}$ with $\nu>1$, then choosing $\phi_1$ and $\psi$ such that $\phi_1^{(j)}(0)\neq 0$ and $\int_{\R} \psi(x_{n-1})dx_{n-1}=1$ it follows from Lebesgue's dominated convergence theorem that
\[
\rho^{2j}e^{-i\eta(\rho^{-1}) \rho T/2}I_{\rho}\to \phi_1^{(j)}(0)/j!\int {\mathcal Y}(x_1;\zeta_0,0)x_1^jdx_1=a_j\neq 0
\]
as $\rho\to\infty$. Note that
\begin{gather*}
|(U(0), D_0v_{\rho}(0))|=\rho^{-2-N/l}|\hat\theta(\rho^5)||I_{\rho}|
=|\hat\theta(\rho^5)|\rho^{-2-N/l-2j}e^{-{\mathsf{Im}}\eta(\rho^{-1})\rho T/2}\\
\times \Big|\rho^{2j}e^{-i\eta(\rho^{-1})\rho T/2}I_{\rho}\Big|\geq (|a_j|/2)|\hat\theta(\rho^5)|\rho^{-2-N/l-2j}e^{-{\mathsf{Im}}\eta(\rho^{-1})\rho T/2},\;\;\rho\geq \rho_0
\end{gather*}
which implies that
\begin{equation}
\label{eq:te:ue}
C\rho^{7+N/l+2j}e^{\rho^{\kappa}+{\mathsf{Im}}\eta(\rho^{-1})\rho T/2}\geq (|a_j|/2)|\hat\theta(\rho^5)|,\quad \rho\geq \rho_0.
\end{equation}
Noting \eqref{eq:stos'} we choose
\[
\phi_1\in \ga_0^{(s_1)}(\R),\;\;\tilde\phi\in \ga_0^{(s_1)}(\R^{n-3}),\;\;\psi\in \ga_0^{(s_2)}(\R),\quad \theta\in \ga_0^{(s_1)}(\R),\;\;\theta\not\in \ga_0^{(5)}(\R)
\]
such that 
\[
\int\tilde\phi(x'') dx''=1,\;\int\psi(x_{n-1})dx_{n-1}=1,\; {\rm supp}\,\phi_1(x_1)\tilde\phi(x'')\theta(x_n)\psi(x_{n-1})\subset L.
\]
It is clear that $\phi_1(x_1)\tilde\phi(x'')\theta(x_n)\psi(\rho^{N/l}x_{n-1})\in \ga_0^{(s),h}(L)$, $\rho\geq 1$ from \eqref{eq:koko}. 
Recalling that ${\mathsf{Im}}\,\eta(0)={\mathsf{Im}}\,\zeta_0<0$ it follows from \eqref{eq:te:ue} that there are $c>0$, $C>0$ such that
\[
Ce^{-c\rho}\geq |\hat\theta(\rho^5)|,\quad \rho\geq \rho_0.
\]
Assuming $\theta$ is even function we conclude that $|\hat\theta(\rho)|\leq Ce^{-c\rho^{1/5}}$, that is $\theta\in \ga_0^{(5)}$ which is a contradiction, thus the proof is completed.
%%%%%%%%%%%%%%

%%%%%%%%%%%%%%
\end{document}